\documentclass[11pt,reqno]{amsart}
\usepackage[T1]{fontenc}
\usepackage{amsmath,amssymb,amsthm,mathtools}
\usepackage[a4paper,left=26mm,right=26mm,top=24.5mm,bottom=24.5mm]{geometry}
\usepackage{xcolor}
\usepackage{needspace}
\usepackage[colorlinks=true,linkcolor=blue!45!black,citecolor=blue!45!black,urlcolor=blue!45!black,pdftitle={Norm descent and Tignol homology over quadric function fields}]{hyperref}
\newtheorem{theorem}{Theorem}[section]
\newtheorem{proposition}[theorem]{Proposition}
\newtheorem{lemma}[theorem]{Lemma}
\newtheorem{corollary}[theorem]{Corollary}
\theoremstyle{remark}
\newtheorem{remark}[theorem]{Remark}
\DeclareMathOperator{\Br}{Br}
\DeclareMathOperator{\Pic}{Pic}
\DeclareMathOperator{\Div}{Div}
\DeclareMathOperator{\Prin}{Prin}
\DeclareMathOperator{\Res}{Res}
\DeclareMathOperator{\Gal}{Gal}

\DeclareMathOperator{\cor}{cor}
\DeclareMathOperator{\disc}{disc}
\DeclareMathOperator{\SB}{SB}
\newcommand{\F}{\mathbb F}
\newcommand{\Z}{\mathbb Z}
\newcommand{\N}{\mathrm N}
\newcommand{\A}{\mathcal A}
\newcommand{\K}{\mathcal K}
\newcommand{\Lam}{\Lambda}
\DeclareMathOperator{\divis}{div}
\DeclareMathOperator{\lc}{lc}
\allowdisplaybreaks[2]
\title[Norm descent and Tignol homology]{Norm descent and Tignol homology over quadric function fields}
\author{Shisong Xu}
\address{Department of Mathematics, Nanjing University, 22 Hankou Road, Nanjing 210093, Jiangsu, People's Republic of China}
\email{shsxu@smail.nju.edu.cn}
\hypersetup{pdfauthor={Shisong Xu}}
\date{}
\begin{document}
\begin{abstract}
Let $L/F$ be a triquadratic extension, with $\operatorname{char}F\ne2$. We describe the kernel of restriction on the first homology of Tignol's complex to the function field of a quadric surface. Its dimension over $\F_2$ is at most two for nonsquare discriminant and six for square discriminant. Both bounds are sharp. For square discriminant we give an optimal bound for each possible collection of splitting subfields of the associated quaternion algebra. Every such collection is realised with equality in the bound over a finitely generated field of transcendence degree at most two over $\mathbb Q$. The kernel for square discriminant is generated by those for associated forms of nonsquare discriminant. For quadratic forms of dimension at least five, descent of principal divisors proves injectivity, answering a question of Sivatski.
\end{abstract}
\maketitle

\section{Introduction}
An element may be a norm from every quadratic subextension of a triquadratic extension without being a norm from the full extension modulo squares. The group of such obstructions can have arbitrarily large dimension over $\F_2$. How many independent obstructions can disappear over the function field of a single quadric? We prove uniform bounds: zero for forms of dimension at least five, and two or six for four-dimensional forms of nonsquare or square discriminant, respectively. Our main result refines the latter bound in terms of the splitting subfields of the associated quaternion algebra and proves its optimality for every possible collection of such subfields, over fields of transcendence degree at most two over $\mathbb Q$.

Throughout, $F$ is a field of characteristic different from two, and quadratic forms are nondegenerate. For a quadratic form $q$ of dimension at least three, write $X_q$ for its projective quadric and $F(q)$ for the function field of $X_q$. Let $L/F$ be a triquadratic extension, and choose $a,b,e\in F^\times$ such that
\[
 L=F(\sqrt a,\sqrt b,\sqrt e),\qquad [L:F]=8.
\]
Let $N_1(L/F)$ denote the first homology of Tignol's complex \cite{Tignol}. Shapiro, Tignol and Wadsworth \cite[Corollaries 1.13, 1.15 and Section 5]{STW} described it using quadratic norm groups and constructed nonzero examples. Sivatski \cite[Section 5]{Siv18} proved that its dimension can be arbitrarily large. Related norm quotients arise in the study of transfer ideals \cite{LW} and $R$-equivalence \cite{Gille,AB}.

The symmetric norm description of $N_1(L/F)$ is due to Tignol and was proved by Gille \cite[Proposition 3]{Gille} using \cite[Corollary 1.15]{STW}; see also \cite[Section 1]{Siv18}. Let $V_L\subset F^\times/F^{\times2}$ be the space of square classes defining $L/F$, written additively. For $d\in V_L$ choose a representative $\widetilde d\in F^\times$, with $\widetilde 0=1$. For a subspace $W\subseteq V_L$ write $F_W$ for the corresponding intermediate field. Put
\begin{equation}\label{eq:normgroups}
 \A(L/F)=\frac{F^\times}{F^{\times2}\N_{L/F}(L^\times)},\qquad
 \Lam(L/F)=\frac{\bigcap_{d\in V_L} \N_d(F)}{F^{\times2}\N_{L/F}(L^\times)},
\end{equation}
where $\N_d(F)=\N_{F(\sqrt{\widetilde d})/F}(F(\sqrt{\widetilde d})^\times)$, with $\N_0(F)=F^\times$. By bilinearity of quaternion symbols, it suffices to take the intersection over a basis of $V_L$. The norm description identifies $N_1(L/F)$ with $\Lam(L/F)$, compatibly with scalar extension preserving the triquadratic degree. We write the norm quotients multiplicatively, using $[c]$ for the class of $c\in F^\times$. Their trivial subgroups are denoted by $0$ when viewed as vector spaces over $\F_2$. Quaternion symbols and classes in the Brauer group are written additively; an algebra $Q$ also denotes its Brauer class when no confusion is possible. We write $\SB(Q)$ for its Severi--Brauer variety.

The problem also has a formulation in terms of Brauer groups. Sivatski \cite[Theorem 3.2]{Siv27} identifies $N_1(L/F)$ with
\[
 B(F;a,b,e)=\frac{{}_2\Br(F(\sqrt a,\sqrt e)/F)\cap {}_2\Br(F(\sqrt b,\sqrt e)/F)}{H(F;a,b,e)},
\]
where, for $t\in F^\times$, we put $R_F(t)=\Br(F(\sqrt t)/F)$ and
\[
 H(F;a,b,e)=R_F(e)+(R_F(a)\cap R_F(b))+(R_F(a)\cap R_F(be)).
\]
Here ${}_2\Br(E/F)$ denotes the 2-torsion in the kernel of restriction from $F$ to $E$. The group $B(F;a,b,e)$ is nonzero precisely when there exists a four-dimensional $(a,b)$-nonstandard form of discriminant $e$, in the terminology of \cite{Siv27}. If a Brauer class in the numerator restricts over $K=F(\sqrt e)$ to $(a,u)=(b,u)$, its image in $\Lam(L/F)$ is $[\N_{K/F}(u)]$ (\cite[Remark 3.4]{Siv27}). This isomorphism commutes with scalar extension preserving the triquadratic degree.

Sivatski asks whether
\begin{equation}\label{eq:question}
 B(F;a,b,e)\longrightarrow B(F(q);a,b,e)
\end{equation}
is injective when $\dim q\geq5$ (\cite[Question 3.12]{Siv27}). If the square classes of $a,b,e$ are dependent, both groups are trivial by \cite[Theorem 3.2]{Siv27}. The analogous assertion for four-dimensional forms fails in general. In characteristic zero the assertion in higher dimension follows from the vanishing of the elementary obstruction for smooth projective hypersurfaces of dimension at least three and the resulting injectivity on $H^1$ of tori (\cite[Lemmas 2.1(iii) and 2.2(ix)]{BCS}, applied to the torus in \eqref{eq:torus}). Our proof uses only divisors over a finite Galois extension and applies in every characteristic different from two.

\begin{theorem}\label{thm:high}
If $\dim q\geq 5$, the map \eqref{eq:question} is injective. More generally, for any finite Galois extension $L/F$ and any positive integer $m$, the natural map
\[
 \frac{F^\times}{F^{\times m}\N_{L/F}(L^\times)}\longrightarrow
 \frac{F(q)^\times}{F(q)^{\times m}\N_{L(q)/F(q)}(L(q)^\times)}
\]
is injective.
\end{theorem}

In geometric dimension at least three the hyperplane class generates the geometric Picard group of a smooth quadric. On a quadric surface the two geometric ruling classes give additional norm relations. For a four-dimensional form of nonsquare discriminant $d$, the even Clifford algebra $Q=C_0(q)$ is a quaternion algebra over $K=F(\sqrt d)$. If $K\subset L$, there are precisely three intermediate fields $M$ of $L/F$ with $[M:K]=2$. For an extension $E/F$ linearly disjoint from $L/F$ define
\[
 \K(L/F;E)=\ker\bigl(\Lam(L/F)\longrightarrow\Lam(LE/E)\bigr).
\]
\begin{theorem}\label{thm:surface-intro}
Let $q$ be four-dimensional with nonsquare discriminant $d$. If $K\not\subset L$ or $q$ is isotropic over $F$, then $\K(L/F;F(q))=0$. Suppose otherwise, and list the fields $M_i=K(\sqrt{r_i})\subset L$, with $r_i\in F^\times$, of degree two over $K$ which split $Q$. Write $s$ for their number.
\begin{enumerate}
\item If $s\leq1$, then $\K(L/F;F(q))=0$.
\item If $s\geq2$, fix $M_1$ and, for $2\leq i\leq s$, choose $u_{1i}\in K^\times$ such that
\[
 Q=(r_1,u_{1i})=(r_i,u_{1i}),\qquad
 u_{1i}=x_{1i}^2-r_1r_i y_{1i}^2,\quad x_{1i},y_{1i}\in K.
\]
Then
\[
 \K(L/F;F(q))=\big\langle[\N_{K/F}(u_{1i})]:2\leq i\leq s\big\rangle\subset\Lam(L/F).
\]
In particular, $\dim_{\F_2}\K(L/F;F(q))\leq s-1\leq2$.
\end{enumerate}
\end{theorem}
The common slots exist by the argument in Section~\ref{sec:surface}; their norm classes need not be independent.

For square discriminant the quadric is stably birational to the conic of a quaternion algebra $Q$ over $F$. Its restriction kernel is governed by the splitting subfields of $Q$. Let
\[
 S(Q,L)=\{d\in V_L\setminus\{0\}:Q_{F_{\langle d\rangle}}=0\},
\]
and let $\mathcal T$ be the set of planes (two-dimensional subspaces) $W\subset V_L$ for which $Q_{F_W}$ splits. The following theorem gives the optimal bound for each collection of splitting subspaces closed under passing to larger subspaces.
\begin{theorem}\label{thm:pattern}
Let $Q$ be a quaternion division algebra over $F$, and let $q$ be its reduced norm form. Put $S=S(Q,L)$ and $\mathcal T_0=\{W\in\mathcal T:W\cap S=\varnothing\}$. Then
\begin{equation}\label{eq:pattern-bound}
 \dim_{\F_2}\K(L/F;F(q))\le b(S,\mathcal T):=
 \begin{cases}
 \max\{|\mathcal T|-1,0\},&S=\varnothing,\\
 |\mathcal T_0|,&S\ne\varnothing.
 \end{cases}
\end{equation}
Conversely, let $V=\F_2^3$ and let $\mathcal S$ be a collection of nonzero subspaces of $V$ closed under passing to larger subspaces (possibly empty). There exist a finitely generated field $F/\mathbb Q$ of transcendence degree at most two, a triquadratic extension $L/F$ with a specified identification $V_L\simeq V$, and a quaternion division algebra $Q/F$ such that
\[
 Q_{F_W}=0\ \Longleftrightarrow\ W\in\mathcal S
 \qquad(0\ne W\subseteq V),
\]
and equality holds in \eqref{eq:pattern-bound} for the prescribed $S$ and $\mathcal T$. If $b(S,\mathcal T)=0$, one may take $F$ to be a number field; otherwise the construction has transcendence degree two.
\end{theorem}
For fixed $L/F$ and $Q$, further arithmetic relations may make \eqref{eq:pattern-bound} strict; see Section~\ref{sec:rank}. Maximising over splitting planes gives the following bounds. Put $r=\dim_{\F_2}\langle S(Q,L)\rangle$.
\begin{corollary}\label{thm:square-intro}
Let $q$ be four-dimensional of square discriminant, similar to the reduced norm form of a quaternion division algebra $Q$. If $S(Q,L)$ is empty, then $\dim_{\F_2}\K(L/F;F(q))\leq6$. Otherwise,
\[
 \dim_{\F_2}\K(L/F;F(q))\leq
 \begin{cases}
 2^{3-r},&\text{if there is }\ell\in V_L^*\text{ with }\ell(d)=1\text{ for every }d\in S(Q,L),\\
 0,&\text{otherwise}.
 \end{cases}
\]
The bound six and each of the positive bounds $2^{3-r}$ for $r=1,2,3$ are attained over finitely generated fields over $\mathbb Q$ of transcendence degree at most two. In particular, the unrestricted optimal bound for square discriminant is six, and the optimal bound in the presence of a quadratic splitting subfield is four.
\end{corollary}

\begin{corollary}\label{thm:sharp-intro}
There exist a finitely generated field $F/\mathbb Q$, a triquadratic extension $L/F$ and an anisotropic four-dimensional form $q$ of nonsquare discriminant for which
\[
 \K(L/F;F(q))\simeq (\Z/2\Z)^2.
\]
One may take $\operatorname{tr.deg}_{\mathbb Q} F=2$.
\end{corollary}
Theorem~\ref{thm:compatibility} expresses the kernel for square discriminant as a sum of seven kernels for forms of nonsquare discriminant associated with the same quaternion algebra.

The finite divisor presentation of Section~\ref{sec:descent} yields the comparison formulas in Sections~\ref{sec:surface} and~\ref{sec:conic}; Section~\ref{sec:rank} gives a criterion for relations among the resulting classes. The sharpness proof in Section~\ref{sec:sharp} starts from the nonzero class of Gille \cite[Remark 4]{Gille}. A quaternion algebra over $\mathbb Q(\sqrt{-3})(x)$ singles out one nonzero comparison, and specialisation along a curve then produces independent comparisons while preserving the prescribed splitting subfields. This differs from the constructions of large homology groups in \cite[Proposition 5.3 and Theorem 5.5]{Siv18} and of nonzero homology in \cite[Theorem 3.6]{Siv27}: here the object being prescribed is the kernel of restriction to a single quadric function field.

\Needspace{16\baselineskip}
\section*{AI Declaration}
OpenAI's ChatGPT (GPT-6) was used for exploratory mathematical work, including calculations of norm comparison identities, polynomial identities for the explicit quaternion example in Proposition~\ref{prop:onehot}, and computations of incidences and splitting patterns in $\F_2^3$. It also assisted in refining the English text. It played no role in generating the core ideas and proofs of the manuscript.

\section*{Computational checks}
Python and SymPy were used for symbolic checks of the polynomial and norm identities, enumeration of the 459 collections of nonzero subspaces of $\F_2^3$ closed under passing to larger subspaces, and verification of the finite incidence relations used in the compatibility theorem. The proofs are theoretical and do not rely on these computational checks.

\section{Descent and a finite divisor presentation}\label{sec:descent}
Let $L/F$ be finite Galois with group $G$, and let $X$ be smooth, projective and geometrically integral over $F$. For an intermediate field $M$ write $X_M=X\times_F M$. The extension $L(X)/F(X)$ is Galois with group $G$. Norm on rational functions is compatible with pushforward of divisors. Write $\N=\N_{L(X)/F(X)}$ in computations on $X$, and use additive notation for divisors.

\begin{proposition}\label{prop:pic}
Suppose that Cartier divisors $H_1,\ldots,H_r$ on $X$ pull back to a $\Z$-basis of $\Pic(X_M)$ for every intermediate field $M$ of $L/F$. Then the inclusion $L^\times\hookrightarrow L(X)^\times$ admits a $G$-equivariant retraction. In particular, for every positive integer $m$,
\[
 F^\times\cap\big(F(X)^{\times m}\N_{L(X)/F(X)}(L(X)^\times)\big)
 =F^{\times m}\N_{L/F}(L^\times).
\]
\end{proposition}
\begin{proof}
For each $G$-orbit of prime divisors of $X_L$, choose a representative $D$ with stabiliser $G_D$ and fixed field $M_D=L^{G_D}$. By Galois descent, $D$ is the pullback of a prime divisor on $X_{M_D}$, denoted by the same letter. By the Picard group hypothesis choose integers $n_{D,j}$ and $g_D\in M_D(X)^\times$ with
\[
 \divis(g_D)=D-\sum_{j=1}^r n_{D,j}H_{j,M_D}.
\]
Extend the assignment to conjugates of $D$ by equivariance; it is well defined since $g_D$ is fixed by $G_D$. Let $s:\Div(X_L)\to L(X)^\times$ be the $G$-equivariant homomorphism which sends each prime divisor to the chosen function. If $D_0=\sum_D a_DD$ is principal, the independence of the classes of the $H_{j,L}$ gives $\sum_D a_D n_{D,j}=0$ for every $j$. Hence $\divis(s(D_0))=D_0$. It follows that
\[
 \rho:L(X)^\times\longrightarrow L^\times,\qquad
 f\longmapsto f/s(\divis f)
\]
is a $G$-equivariant retraction: the quotient has zero divisor and $X_L$ is proper and geometrically integral. On invariants $\rho$ takes values in $F^\times$ and commutes with the norm. Applying $\rho$ to $c=h^m\N_{L(X)/F(X)}(f)$ for $c\in F^\times$ proves the nontrivial inclusion.
\end{proof}

\begin{remark}
This is a finite Galois version of the retraction underlying the elementary obstruction in \cite[Section 2]{BCS}. For any $F$-torus $T$ split by $L$, tensoring the retraction with its cocharacter lattice gives a $G$-equivariant retraction $T(L(X))\to T(L)$. Hilbert 90 and inflation--restriction therefore give injectivity of $H^1(F,T)\to H^1(F(X),T)$ under the hypotheses of the proposition.
\end{remark}

\begin{proof}[Proof of Theorem~\ref{thm:high}]
For a smooth quadric of dimension at least three the geometric Picard group is generated by the hyperplane class. This class is defined over the ground field, so it is a basis over every intermediate field as well; see \cite{EKM}. Proposition~\ref{prop:pic} gives the claimed injectivity on norm quotients. Taking $m=2$ and using \eqref{eq:normgroups} proves the assertion for $B(F;a,b,e)$.
\end{proof}

\begin{remark}
The argument also gives injectivity on the norm quotient for a rational function field in any finite number of variables, by taking $X$ to be projective space. Consequently the kernel considered in this paper is unchanged by stable birational equivalence of smooth projective geometrically integral varieties.
\end{remark}

For any finite Galois extension $L/F$ with group $G$, let $T$ be the $F$-torus defined by the exact sequence
\begin{equation}\label{eq:torus}
 1\longrightarrow T\longrightarrow \mathbf G_m\times\Res_{L/F}\mathbf G_m
 \xrightarrow{\ (z,w)\mapsto z^2\N_{L/F}(w)\ }\mathbf G_m\longrightarrow1.
\end{equation}
The map on character lattices is primitive, so $T$ is a torus. Hilbert 90 identifies $H^1(F,T)$ with $\A(L/F)$, compatibly with extension to $F(X)$. Let $M=X_*(T)$ be its cocharacter lattice, $\pi:X_L\to X$ the projection, and $\mathcal P_E=\Prin(X_E)$ the group of principal divisors over an intermediate field $E$; write $\mathcal P_F=\Prin(X)$. Here $\A(L/F)$ is defined by \eqref{eq:normgroups} for any finite Galois extension.

\begin{proposition}\label{prop:canonical}
There is a canonical exact sequence
\begin{equation}\label{eq:canonical}
 T(F(X))\xrightarrow{\divis} (M\otimes_{\Z}\mathcal P_L)^G
 \xrightarrow{\delta}\A(L/F)\longrightarrow\A(L(X)/F(X)).
\end{equation}
Here the middle group can be identified with
\[
 \mathcal R_X=\{D\in\mathcal P_L:\pi_*D\in2\mathcal P_F\}.
\]
Under this identification the first map sends $(z,f)$ to $\divis(f)$, and the second sends $D=\divis(f)$ to
\[
 \left[\N_{L(X)/F(X)}(f)/h^2\right],
 \qquad \divis(h)=\tfrac12\pi_*D.
\]
In particular the kernel on norm quotients is canonically isomorphic to the cokernel of the first map in \eqref{eq:canonical}. This sequence requires no assumption on the torsion of $\Pic(X)$.
\end{proposition}
\begin{proof}
The torus $T$ splits over $L$. Its cocharacter lattice is the kernel of the surjection
\[
 \Z\oplus\Z[G]\longrightarrow\Z,\qquad (n,(n_g))\longmapsto 2n+\sum_{g\in G}n_g.
\]
The exact sequence of $G$-modules
\[
 0\longrightarrow L^\times\longrightarrow L(X)^\times
 \xrightarrow{\divis}\mathcal P_L\longrightarrow0
\]
remains exact after tensoring with the free lattice $M$. The resulting first two modules identify with $T(L)$ and $T(L(X))$. Taking $G$-cohomology gives \eqref{eq:canonical}, since $H^1(L,T)=H^1(L(X),T)=0$.

The lattice description of $M$ identifies its tensor with $\mathcal P_L$ with tuples $(D_0,(D_g))$ satisfying $2D_0+\sum_gD_g=0$. A $G$-invariant tuple is specified by its component $D$ at the identity of $G$, with $D_g=gD$ and $D_0\in\mathcal P_L^G$. Galois descent for divisors and Hilbert 90 identify $\mathcal P_L^G$ with the pullback of $\mathcal P_F$. The identity $\sum_g gD=\pi^*\pi_*D$ therefore gives the claimed identification with $\mathcal R_X$. Lifting the tuple to rational functions in the defining sequence of $T$ gives the constant $\N(f)/h^2$ (its inverse represents the same class).
\end{proof}

To obtain a finite presentation, assume that $G$ is abelian and that $\Pic(X)$ is free of 2-torsion. Let $M_1,\ldots,M_n$ be intermediate fields (repetitions allowed) and $D_i\in\Div(X_{M_i})$. Assume that for every intermediate $M$ the classes of pullbacks of the $G$-conjugates of $D_i$ with $M_i\subseteq M$ generate $\Pic(X_M)$. The fields $F$ and $L$ are allowed among the $M_i$.

Introduce the permutation lattice
\[
 E=\bigoplus_{i=1}^n\Z[G/\Gal(L/M_i)]
\]
and let $\mathcal D:E\to\Div(X_L)$ send each basis element to the corresponding conjugate of the pullback of $D_i$. Write $p:E\to\Pic(X_L)$ for its composite with the divisor class map. Define
\[
 \nu:E\longrightarrow(\Z/2\Z)^n,\qquad
 \nu(z)_i=[L:M_i]\sum_{\sigma\in G/\Gal(L/M_i)}z_{i,\sigma}\pmod 2,
\]
and let $W=\ker p\cap\ker\nu$. For $z\in W$, choose $f_z\in L(X)^\times$ with $\divis(f_z)=\mathcal D(z)$. Since the pushforward of $\mathcal D(z)$ is even, write it as $2D_z$ on $X$. The class of $2D_z$ vanishes, so the 2-torsion hypothesis gives $D_z=\divis(h_z)$ for some $h_z\in F(X)^\times$. Then
\begin{equation}\label{eq:theta}
 \theta(z)=\left[\frac{\N_{L(X)/F(X)}f_z}{h_z^2}\right]\in\A(L/F).
\end{equation}
The quotient in brackets is constant. Changing the choices modifies it by a square and a norm from $L$. Thus $\theta$ is a well-defined homomorphism.

\begin{proposition}\label{prop:finite}
Under the hypotheses above,
\[
 \ker\bigl(\A(L/F)\longrightarrow\A(L(X)/F(X))\bigr)=\theta(W).
\]
If $\Pic(X_L)$ is free abelian, $\theta$ induces a surjection from the finite-dimensional $\F_2$-vector space
\[
 \ker\bigl(p\bmod2:E/2E\to\Pic(X_L)/2\Pic(X_L)\bigr)
 \cap\ker(\nu\bmod2).
\]
\end{proposition}
\begin{proof}
The inclusion of $\theta(W)$ in the kernel is immediate. Conversely, suppose $c\in F^\times$ has the form
\[
 c=\N_{L(X)/F(X)}(f)h^2.
\]
For each orbit of prime divisors in $\divis(f)$ choose $D$ and $M_D$ as in Proposition~\ref{prop:pic}. Express $[D]\in\Pic(X_{M_D})$ as a linear combination of reference classes defined over $M_D$. The difference of the corresponding divisors is $\divis(g_D)$ for some $g_D\in M_D(X)^\times$. Extend this choice equivariantly. Dividing $f$ by the corresponding product of these functions gives $f_0$ with $\divis(f_0)=\mathcal D(z)$ for some $z\in\ker p$.

For the finite constant extension $X_L\to X$, the prime divisors over a fixed prime divisor of $X$ form one $G$-orbit. Each pushes forward with multiplicity $|G_D|=[L:M_D]$: the extension is finite \'etale, and this is its residue field degree at the chosen prime. Since $\divis(\N(f))=-2\divis(h)$, the sum of the exponents in each orbit, multiplied by $[L:M_D]$, is even. More explicitly, if the coefficients of the conjugates of $D$ are $a_\sigma$, the norm of their removed functions is
\[
 \N_{M_D(X)/F(X)}(g_D)^{[L:M_D]\sum_\sigma a_\sigma}.
\]
It is a square. If a reference divisor over $M_i$ occurs with coefficient $b_i$ in the expression of $[D]$, its contribution to $\nu(z)_i$ is $[L:M_i]b_i\sum_\sigma a_\sigma$ modulo two (sum these contributions if several conjugates occur). This is zero because $M_i\subseteq M_D$. Summing over the orbits proves the parity assertion. Thus $z\in W$, and $[c]=\theta(z)$ by \eqref{eq:theta}.

Since $p$ is surjective and $\Pic(X_L)$ is free abelian, $\ker p$ reduces onto $\ker(p\bmod2)$. Thus every vector in the stated intersection lifts to $W$. If $z=2w\in W$, then $p(w)=0$. A rational function for $\mathcal D(z)$ can be chosen to be the square of one for $\mathcal D(w)$; its norm is a square, so $\theta(z)=1$. Hence $\theta$ factors through the displayed finite-dimensional space.
\end{proof}

\begin{lemma}\label{lem:relations}
Suppose $z\in W$ and $\mathcal D(z)$ is invariant under a nontrivial subgroup $J$ of $G$, with $|J|$ even. Then $\theta(z)=1$. Also $\theta(\sigma z)=\theta(z)$ for every $\sigma\in G$.
\end{lemma}
\begin{proof}
The invariant divisor descends to $X_{L^J}$ and is principal there: the restriction $\Pic(X_{L^J})\to\Pic(X_L)$ is injective. This follows from Hilbert 90 and $H^0(X_L,\mathcal O_{X_L}^\times)=L^\times$. The function $f_z$ may therefore be chosen over $L^J$, so its norm to $F(X)$ is a square. Its quotient by $h_z^2$ is a square in $F(X)$ and hence in $F$. The second assertion follows from invariance of norm under conjugation.
\end{proof}

\section{Quadric surfaces of nonsquare discriminant}\label{sec:surface}
Every four-dimensional form of discriminant $d$ is similar to a form
\begin{equation}\label{eq:normalform}
 \langle d,-\alpha,-\beta,\alpha\beta\rangle,
 \qquad \alpha,\beta\in F^\times.
\end{equation}
For $\langle a_1,a_2,a_3,a_4\rangle$, take $\lambda=a_4/(a_2a_3)$, $\alpha=-\lambda a_2$ and $\beta=-\lambda a_3$. Then $\lambda a_4=\alpha\beta$ and $\lambda a_1$ has square class $d$. For nonsquare $d$ its even Clifford algebra has centre $K=F(\sqrt d)$ and is $(\alpha,\beta)_K$: if $e_2,e_3,e_4$ are the Clifford generators for the last three entries, then $e_2e_4/\alpha$ and $e_3e_4/\beta$ anticommute and have squares $\beta$ and $\alpha$, respectively. They generate the algebra over its centre. The even Clifford algebra is invariant under similarity (see \cite{EKM}). The form is isotropic precisely when this algebra splits. To see the nontrivial direction, write $Q_0=(\alpha,\beta)_F$. If $Q_0$ splits, the subspace of elements of trace zero contains a nonzero nilpotent. If $Q_0$ is division and $(Q_0)_K$ splits, $K$ embeds in $Q_0$ (by the maximal subfield criterion for a central simple algebra; see \cite{GS}) and the pure element with square $d$ has reduced norm $-d$. The last three coordinates in \eqref{eq:normalform} are the reduced norm on the subspace of elements of trace zero of $Q_0$, so either case gives an isotropic vector. Conversely, an isotropic four-dimensional form becomes hyperbolic over its discriminant field and hence has split even Clifford algebra there. For square discriminant the same normal form is the reduced norm form of $Q_0$.

We use the common slot lemma of Sivatski \cite[Lemma 2.1]{Siv27}. Suppose a quaternion division algebra $Q/K$ is split by $K(\sqrt r)$ and $K(\sqrt s)$, where $r,s$ have independent square classes. Choose elements of trace zero with squares $r,s$ in $Q$. In the three-dimensional space of elements of trace zero, choose a nonzero element orthogonal to both for the polar form of the reduced norm. It anticommutes with both elements and has nonzero square $u$, since the algebra is division. Thus $Q=(r,u)=(s,u)$, so $(rs,u)=0$ and $u=x^2-rs y^2$ for $x,y\in K$. If $Q$ is split, take $u=1$.

Let $q$ be four-dimensional of nonsquare discriminant $d$, let $K=F(\sqrt d)$, and put $Q=C_0(q)$ over $K$. If $K\not\subset L$, the Picard groups of $X_{q,M}$ for all intermediate $M$ are generated by the hyperplane class. Indeed, the geometric Picard group has basis the two ruling classes, exchanged by the discriminant character. Its invariant subgroup when the discriminant is nonsquare is generated by their sum, which descends. Proposition~\ref{prop:pic} gives injectivity in this case. If $q$ is isotropic, $X_q$ is rational and the result follows in the same way. We henceforth assume $K\subset L$ and $q$ anisotropic.

Over $K$ the form has square discriminant and is similar to the norm form of a quaternion algebra of class $Q$. Accordingly $X_{q,K}$ is the product of the conics of $Q$ and $Q^{\mathrm{op}}$. Over a splitting field, a nonzero matrix of reduced norm zero of size two has rank one, and its image and kernel determine its projective class. This identifies the quadric with the product of the two varieties of one-dimensional reduced ideals; the identification descends by equivariance under change of splitting. The two conics have the same Brauer class. In geometric ruling coordinates $(u,v)$ the Picard groups over the intermediate fields containing $K$ are
\begin{equation}\label{eq:picquadric}
 \Pic(X_{q,M})=\begin{cases}
 \Z^2,&Q_M=0,\\
 \{(u,v)\in\Z^2:u+v\equiv0\pmod2\},&Q_M\ne0.
 \end{cases}
\end{equation}
Indeed the obstruction to descending $\mathcal O(u,v)$ on the product is $(u+v)Q_M$ in $\Br(M)$. This follows from the Picard--Brauer exact sequence and the Brauer obstruction of a ruling (see \cite[Corollary 2.5.9 and Section 6.1]{CTS}). Over intermediate fields not containing $K$, the Picard group is generated by the hyperplane class $H=(1,1)$. By anisotropy, $Q$ is division over $K$. We choose a divisor $U$ on $X_{q,K}$ of class $(2,0)$, pulled back from a degree-two divisor on the first conic. For $\tau\in G\setminus\Gal(L/K)$ let $U'=\tau U$, of class $(0,2)$.

The geometric Picard invariants are generated by the hyperplane class, which descends to $F$. The Picard--Brauer sequence and the injectivity of $\Br(X_q)\to\Br(F(q))$ (\cite[Section 3.5]{CTS}) give injectivity of $\Br(F)\to\Br(F(q))$. In particular a constant $c\in F^\times$ representing an element of the kernel on norm quotients lies in all the quadratic norm groups in \eqref{eq:normgroups}, since $(\widetilde r,c)$ splits over $F(q)$ for every $r\in V_L$. Thus in this section the kernel on $\A$ is exactly the kernel on $\Lam$.

\begin{proof}[Proof of Theorem~\ref{thm:surface-intro}]
We distinguish three cases. All divisors mentioned below are pulled back to $X_{q,L}$ when a relation is formed; reference fields are understood as in Section~\ref{sec:descent}.

\emph{Case 1: $Q_L$ is nonsplit.} The reference divisors $H,U$ over $F,K$, respectively, and their conjugates satisfy the hypotheses of Proposition~\ref{prop:finite} by \eqref{eq:picquadric}. Every reference field is proper in $L$, so $\nu=0$. The relations between the three classes $H,U,U'$ are generated by $U+U'-2H$, which is principal and defined over $K$. Its image under $\theta$ is trivial by Lemma~\ref{lem:relations}.

\emph{Case 2: $Q_L$ splits but no proper intermediate field containing $K$ splits $Q$.} In addition to $H,U,U'$, choose a ruling divisor $P$ of class $(1,0)$ over $L$ and all its conjugates, with $P'=\tau P$ of class $(0,1)$. Adding the two ruling coordinates of a relation shows that the sum of the coefficients of the conjugates of $P$ is even. Thus $W=\ker p$.

Write $J=\Gal(L/K)\simeq(\Z/2\Z)^2$. The relations are generated by
\begin{gather*}
 \sigma P-P,\quad \tau(\sigma P-P)\quad(\sigma\in J),\\
 P+P'-H,\quad 2P-U,\quad 2P'-U',
\end{gather*}
since the two classes of $P,P'$ form a basis. We use the same notation for a relation in the reference lattice and its image in the divisor group. The middle relation $P+P'-H$ is invariant under $\tau$ and hence has trivial image. The map
\[
 J\longrightarrow\A(L/F),\qquad\sigma\longmapsto\theta(\sigma P-P)
\]
is a homomorphism by additivity and conjugation invariance. For every nonidentity $\sigma\in J$, the relation $P+\sigma P-U$ is invariant under $\sigma$ and so has trivial image. Consequently $\theta(\sigma P-P)=\theta(2P-U)$ for all three nonidentity elements of $J$. A homomorphism from $(\Z/2\Z)^2$ to an elementary abelian $2$-group which is constant on its three nonzero elements is zero. All generators therefore have trivial image.

\emph{Case 3: some proper intermediate fields containing $K$ split $Q$.} Let $M_1,\ldots,M_s$ be the splitting fields of degree two over $K$. Choose $P_i$ of class $(1,0)$ over $M_i$, pulled back from a rational point on the first conic over $M_i$. Take their conjugates together with $H,U,U'$. By \eqref{eq:picquadric} these are a valid reference system. Since all reference fields are proper, $W=\ker p$. Let $P_i'=\tau P_i$. Differences of two conjugates of $P_i$ in the same ruling, and the relations $2P_i-U$, $2P_i'-U'$, $P_i+P_i'-H$, are principal over $M_i$ and have trivial image by Lemma~\ref{lem:relations}. Using the classes of $P_1,P_1'$ as a basis, the remaining relations are $P_i-P_1$ and their conjugates. Lemma~\ref{lem:relations} therefore gives the generators $\theta(P_i-P_1)$ for $i\geq2$.

It remains to compute these classes. Fix $i\geq2$, and put $r=r_1$, $s_0=r_i$. The compositum $M_1M_i=L$, so $L=K(\sqrt r,\sqrt{s_0})$. Since $Q$ is split by both quadratic extensions of $K$, the common slot argument above gives
\[
 Q=(r,u)=(s_0,u),\qquad u=x^2-rs_0y^2\ne0,
 \quad x,y\in K.
\]
On the first conic use the equation $X^2-rY^2=uZ^2$ and the rational coordinate $t=(X+\sqrt rY)/Z$ over $L$. The two points at infinity are the points $t=0,\infty$. We may take $P_1$ from $t=\infty$ and $P_i$ from $(X,Y,Z)=(x,\sqrt{s_0}y,1)$; changing a point within its reference field changes $\theta$ by a trivial class. Set $w=x+\sqrt{rs_0}y$ and $w'=x-\sqrt{rs_0}y$, so $ww'=u$. The function $f=t-w$ on the first conic, pulled back to $X_{q,L}$, has divisor $P_i-P_1$. The automorphisms changing the signs of $\sqrt r$ and $\sqrt{s_0}$, respectively, send $t$ to $u/t$ and $t$ to itself; each interchanges $w$ and $w'$. Hence
\begin{align*}
 \N_{L(X_q)/K(X_q)}(f)
 &= (t-w)(t-w')(u/t-w')(u/t-w)\\
 &=u\left(\frac{(t-w)(t-w')}{t}\right)^2\\
 &=u\bigl(2(X/Z-x)\bigr)^2.
\end{align*}
The expression in parentheses is defined over $K$. Norm to $F(X_q)$ therefore gives $\theta(P_i-P_1)=[\N_{K/F}(u)]$. The class is independent of $u,x,y$, since changing either ruling divisor over its field of definition gives a relation with trivial image.
\end{proof}

\begin{corollary}\label{cor:nonstandard}
Suppose $q$ is four-dimensional, of discriminant $e$, and is isotropic over both $F(\sqrt a)$ and $F(\sqrt b)$. If $q$ is $(a,b)$-nonstandard, then $\K(L/F;F(q))\ne0$. If its even Clifford algebra is not split by $F(\sqrt{ab},\sqrt e)$, this kernel has dimension exactly one, generated by the norm class of the obstruction of $q$.
\end{corollary}
\begin{proof}
For a nonstandard form let $Q_0$ be a quaternion algebra over $F$ restricting to $C_0(q)$ over $F(\sqrt e)$. Its image in $B(F;a,b,e)$ is nonzero by the definition in \cite{Siv27}. Over $F(q)$ the form is isotropic, hence standard, so this image belongs to the restriction kernel. The two fields $F(\sqrt a,\sqrt e)$ and $F(\sqrt b,\sqrt e)$ split the even Clifford algebra. If the third field does not, Theorem~\ref{thm:surface-intro} gives the bound one. The norm representative of the nonzero obstruction is the same $[\N_{K/F}(u)]$ as in the proof, by the explicit norm comparison of \cite[Remark 3.4]{Siv27}.
\end{proof}

\begin{proposition}\label{prop:example}
Let $F=\mathbb Q(x)$ and take
\[
 a=x+4,\qquad b=x+1,\qquad e=x,\qquad c=2(x+2).
\]
For the four-dimensional form $q=\langle e,-a,-c,ac\rangle$ and the triquadratic extension $L=F(\sqrt a,\sqrt b,\sqrt e)$ one has
\[
 \K(L/F;F(q))=\langle[-x]\rangle\simeq\Z/2\Z
\]
in the norm description \eqref{eq:normgroups}.
\end{proposition}
\begin{proof}
The square classes of $a,b,e$ are independent by the valuations at their distinct zeros. Write $K=F(s)$, $s^2=x$, and put $Q=(a,c)_K=C_0(q)$. We have
\[
 s=\left(\frac{s+2}{2}\right)^2-a\left(\frac12\right)^2,
 \qquad u=sc=(s^2+s+2)^2-ab.
\]
The first identity gives $(a,s)=0$ over $K$ and the second gives $(ab,u)=0$. Consequently $Q=(a,u)=(b,u)$, so $Q$ splits over $K(\sqrt a)$ and $K(\sqrt b)$. The norm of the common slot is
\[
 \N_{K/F}(u)=-4x(x+2)^2.
\]
Its class in $\Lam(L/F)$ is $[-x]$. This is the classical nonzero class of \cite[Remark 4]{Gille}, based on \cite[Theorem 5.1 and Remark 5.4]{STW}. In fact, the class of $2=((s+2)^2-a)/((s+1)^2-b)$ in the quotient of \cite[Corollary 1.15]{STW} maps to $[\N_{K/F}(4s)]=[-x]$ by \cite[Proposition 3]{Gille}. Nontriviality already holds over $\mathbb Q_2(x)$.

To see that $Q$ does not split over $K(\sqrt{ab})$, use the valuation of $K=\mathbb Q(s)$ defined by $s^2+2$. Its residue field is $\mathbb Q(\sqrt{-2})$. Here $a$ is a unit with residue $2$, while $c$ has valuation one. Thus the residue of $Q$ is the nonsquare class of $2$ in $\mathbb Q(\sqrt{-2})^\times/\mathbb Q(\sqrt{-2})^{\times2}$. On the other hand $ab$ is a unit with square residue $-2$, so the quadratic extension $K(\sqrt{ab})/K$ splits at the completion of this valuation and the nonzero residue persists. Hence the third quadratic field over $K$ does not split $Q$. Theorem~\ref{thm:surface-intro} now gives the assertion.
\end{proof}

\section{Conics and square discriminant}\label{sec:conic}
We keep the triquadratic extension $L/F$ fixed. Let $C=\SB(Q)$ for a quaternion algebra $Q$ over $F$. Write
\[
 J_C=\ker\bigl(\A(L/F)\longrightarrow\A(L(C)/F(C))\bigr).
\]
If $Q$ is split over $F$, the kernel is trivial by rationality. If $Q_L$ is nonsplit, every intermediate Picard group is generated by the hyperplane class of the plane conic $C$, so the kernel is again trivial by Proposition~\ref{prop:pic}.

Suppose that $Q$ is division over $F$ and split over $L$. Let $M_1,\ldots,M_s$ be the minimal intermediate fields which split $Q$. Choose a point $P_i\in C(M_i)$ and consider it as a divisor of degree one on $C_L$. Choose a hyperplane divisor $U$ on $C$ of degree two. If all $M_i$ are proper subfields of $L$, define, for $i\ne j$,
\begin{equation}\label{eq:cij}
 c_{ij}=\left[\frac{\N_{L(C)/F(C)}(f_{ij})}{h_{ij}^2}\right]\in\A(L/F),\qquad
 \divis(f_{ij})=P_i-P_j,
\end{equation}
where $h_{ij}$ is a function over $F$ whose divisor is half the pushforward of $P_i-P_j$. The pushforward is even and $\Pic(C)$ is torsion-free, so $h_{ij}$ exists. Changing the points gives a principal divisor over a proper subfield, whose defining function has square norm from $L$. Changing the functions without changing their divisors multiplies the quotient by a constant norm or square. The same definition applies to any pair of proper splitting subfields, not necessarily minimal. We write $c_{M,N}$ for this comparison class when the fields are $M,N$. With $c_{M,M}=1$, we have $c_{M,N}c_{N,P}=c_{M,P}$ and $c_{M,N}=1$ if $MN\ne L$. These assertions follow by multiplying the functions in \eqref{eq:cij} and, in the latter case, choosing the principal function over $MN$. These comparisons are compatible with extensions of the ground field for which $L$ remains a field.

\begin{theorem}\label{thm:conic}
If $Q$ is division over $F$, the following statements hold.
\begin{enumerate}
\item If $Q_L$ is nonsplit, or $L$ itself is the only minimal splitting subfield, then $J_C=0$.
\item Otherwise, with the notation above, $J_C=\langle c_{i1}:2\leq i\leq s\rangle$. In particular $\dim_{\F_2}J_C\leq s-1$. Furthermore, $c_{ij}=1$ whenever $M_iM_j\ne L$.
\end{enumerate}
\end{theorem}
\begin{proof}
The case in which $Q_L$ is nonsplit has already been proved. If the only splitting field is $L$, choose one point $P\in C(L)$, and use $U$ and the conjugates of $P$ as reference divisors. Over any proper intermediate field the Picard group of the conic is generated by $U$, whereas the point class generates the Picard group over $L$. In any relation the sum of the coefficients of the conjugates of $P$ is even, so $W=\ker p$. This relation group is generated by $\sigma P-P$ and $2P-U$ for $\sigma\in G$. For every nonidentity $\sigma\in G$, $P+\sigma P-U$ is an invariant principal divisor under $\sigma$ and has trivial image by Lemma~\ref{lem:relations}. The homomorphism $\sigma\mapsto\theta(\sigma P-P)$ is therefore constant on all nonidentity elements, with value $\theta(2P-U)$. Since $G$ has two independent elements of order two, this value is trivial. Proposition~\ref{prop:finite} yields $J_C=0$.

If there is a proper splitting field, all minimal splitting fields are proper. Use $U$ and all conjugates of the $P_i$ as reference divisors. The hypotheses of Proposition~\ref{prop:finite} hold: over a nonsplitting intermediate field $U$ is a generator, and over any splitting intermediate field a point from a minimal splitting subfield gives a generator. All reference fields are proper, so $W=\ker p$. Differences of conjugates in a single reference field, and $2P_i-U$, have trivial image. Subtracting multiples of the class of $P_1$ shows that the remaining relations are generated by $P_i-P_1$ and their conjugates. Their images are exactly the classes in \eqref{eq:cij}. When $M_iM_j$ is proper, the relation $P_i-P_j$ is principal over the compositum and its norm is a square. This proves the assertions.
\end{proof}

Consider the graph with vertices the proper minimal splitting fields, and an edge between two fields when their compositum is proper in $L$. If there are $b>0$ connected components, the theorem and the comparison identity give a surjection from the reduced degree-zero homology of this graph over $\F_2$ onto $J_C$. Indeed, differences within each component have trivial image and differences from one fixed vertex generate $J_C$. Consequently $\dim_{\F_2}J_C\le b-1$.

We next show $J_C\subset\Lam(L/F)$. Let $M,N$ be distinct biquadratic splitting subfields, put $K_0=M\cap N$, and write $M=K_0(\sqrt r)$, $N=K_0(\sqrt s)$, with $r,s\in F^\times$. The common slot argument in Section~\ref{sec:surface} gives
\[
 Q_{K_0}=(r,u)=(s,u),\qquad u=x^2-rsy^2\in K_0^\times,
\]
where $x,y\in K_0$ (if $Q_{K_0}=0$, take $u=x=1$ and $y=0$). The calculation on the conic in Section~\ref{sec:surface}, now applied to $C_{K_0}$, gives
\begin{equation}\label{eq:pairnorm}
 c_{M,N}=[\N_{K_0/F}(u)].
\end{equation}
This class belongs to $\Lam(L/F)$: it is a norm from $K_0$, and the projection formula gives $(r,\N_{K_0/F}u)=\cor_{K_0/F}(Q_{K_0})=2Q=0$, and likewise for $s$. A radicand for $K_0$ together with $r,s$ gives a basis of $V_L$. Each proper splitting subfield is contained in a biquadratic splitting subfield, and its comparison with that larger field is trivial. Therefore all comparisons in Theorem~\ref{thm:conic} reduce to \eqref{eq:pairnorm} or to the trivial class.

If $q$ has dimension four and square discriminant, it is similar to the reduced norm form of a quaternion algebra $Q$ over $F$. Its quadric is the product of the conics of $Q$ and $Q^{\mathrm{op}}$; hence $F(q)$ is rational of transcendence degree one over $F(C)$ for $C=\SB(Q)$. By rational invariance and the inclusion just proved, $\K(L/F;F(q))=J_C$.
\begin{corollary}\label{cor:square}
For a four-dimensional form of square discriminant and a triquadratic extension $L/F$,
\[
 \K(L/F;F(q))=J_C.
\]
In particular the classes in \eqref{eq:cij} generate this kernel whenever the proper minimal splitting fields exist. If $Q$ is division and $t$ denotes the number of minimal biquadratic splitting subfields, then:
\begin{enumerate}
\item if there is a quadratic splitting subfield, the kernel has dimension at most $t\leq4$;
\item if there is no quadratic splitting subfield but $t>0$, it has dimension at most $t-1\leq6$;
\item if there is no proper splitting subfield, it is zero.
\end{enumerate}
\end{corollary}
\begin{proof}
Only the bounds require explanation. When there are quadratic splitting fields, any two of them have proper compositum, so their comparison class is trivial by Theorem~\ref{thm:conic}. Choose one of them as the reference field. Only the minimal biquadratic splitting fields then contribute generators, and at most four of the seven biquadratic subfields avoid a given quadratic subfield. In the absence of quadratic splitting fields the only proper minimal ones are biquadratic; choosing one as reference gives at most $t-1$ generators. The remaining case is part of Theorem~\ref{thm:conic}.
\end{proof}

\begin{corollary}\label{cor:refined}
With the notation of Corollary~\ref{thm:square-intro}, suppose $S(Q,L)\ne\varnothing$. The dimension of $J_C$ is at most the number of planes in $V_L$ disjoint from $S(Q,L)$. This number is $2^{3-r}$ if the linear equations $\ell(d)=1$ for $d\in S(Q,L)$ are consistent and zero otherwise.
\end{corollary}
\begin{proof}
All quadratic splitting fields have trivial comparisons with one another, since the compositum of any two is proper in $L$. A biquadratic splitting field containing one of them also has trivial comparison with each of them. Thus only planes avoiding all quadratic splitting classes can contribute nonzero generators to the formula in Theorem~\ref{thm:conic}. Each plane is the kernel of a unique nonzero functional on $V_L$. Avoiding $S(Q,L)$ amounts precisely to the displayed linear equations. Since the set of constraints is nonempty, a solution is automatically nonzero; when consistent the solution space is affine of dimension $3-r$. This proves the assertion and the upper bounds in Corollary~\ref{thm:square-intro}.
\end{proof}

Assume $Q$ is division; otherwise the following kernels vanish by rationality. Fix a presentation $Q=(\alpha,\beta)_F$. For each $d\in V_L$ put
\[
 q_d=\langle \widetilde d,-\alpha,-\beta,\alpha\beta\rangle,
 \qquad I_d=\K(L/F;F(q_d)).
\]
Thus $I_0=J_C$. For $d\ne0$, the even Clifford algebra of $q_d$ is $Q_{F(\sqrt{\widetilde d})}$. Write $c_{W,W'}=c_{F_W,F_{W'}}$ for splitting subspaces, and let $\mathcal P$ be the set of all seven planes in $V_L$.

\begin{theorem}\label{thm:compatibility}
For $d\ne0$,
\[
 I_d=\langle c_{W,W'}: W,W'\in\mathcal T,\ d\in W\cap W'\rangle,
 \qquad J_C=\sum_{d\ne0}I_d.
\]
If $\dim_{\F_2}J_C=6$, all the $I_d$ with $d\ne0$ have dimension two. For distinct nonzero $d,e$ one has
\[
 I_d\cap I_e=0,\qquad J_C=I_d\oplus I_e\oplus I_{d+e}.
\]
For linearly independent $d,e,f$, the sum $I_d+I_e+I_f$ has dimension five.
\end{theorem}
\begin{proof}
When $Q_{F(\sqrt{\widetilde d})}$ is nonsplit, the first formula follows from Theorem~\ref{thm:surface-intro} and \eqref{eq:pairnorm}, since the quadratic splitting fields over $F(\sqrt{\widetilde d})$ are exactly the fields of the splitting planes through $d$. When it is split, $q_d$ is isotropic, and every comparison between planes through $d$ is trivial by comparison with $F(\sqrt{\widetilde d})$. Every pair of distinct planes meets in a unique line, and the comparisons between splitting planes generate $J_C$ by Theorem~\ref{thm:conic} and the reduction to planes preceding Corollary~\ref{cor:square}. This proves the sum formula.

In dimension six all seven planes split $Q$, and there is no quadratic splitting subfield, by Corollary~\ref{cor:square}. The comparison identities give a surjection from the augmentation space
\[
 A_0(\mathcal P)=\{(x_W)\in\F_2^{\mathcal P}:\textstyle\sum_Wx_W=0\}
 \longrightarrow J_C,\qquad e_W+e_{W'}\longmapsto c_{W,W'},
\]
which is an isomorphism by dimension. The first formula identifies $I_d$ with the augmentation space supported on the three planes through $d$. In particular, these subspaces are independent of the presentation of $Q$. Two such sets meet in one plane, so their augmentation spaces intersect trivially. The sum of the augmentation spaces of a family of these sets is the augmentation space on their union: each pair of sets meets, and differences of basis vectors along the resulting connected graph generate all differences in the union. For $d,e,d+e$ the union is all seven planes; for independent $d,e,f$ it is six planes. The stated dimensions and direct sum follow.
\end{proof}

\begin{corollary}\label{cor:product}
Suppose all seven biquadratic subfields of $L/F$ split $Q$. For a subset $D\subset V_L\setminus\{0\}$ let $E_D=F(\prod_{d\in D}X_{q_d})$, taking $E_\varnothing=F$. Then the restriction on norm quotients has kernel $\sum_{d\in D}I_d$ and induces an exact sequence
\[
 0\longrightarrow\sum_{d\in D}I_d\longrightarrow J_C
 \longrightarrow J_{C_{E_D}}\longrightarrow0,
\]
where the last kernel is formed with $LE_D/E_D$. In particular, if $\dim_{\F_2} J_C=6$ and $D\ne\varnothing$, then
\[
 \dim_{\F_2}J_{C_{E_D}}=7-\#\{W\in\mathcal P:W\cap D\ne\varnothing\}.
\]
\end{corollary}
\begin{proof}
Adjoin the function fields of the quadrics successively. All seven planes remain splitting planes, and the radicands remain independent by geometric integrality. At each stage Theorem~\ref{thm:compatibility} identifies the next kernel with the image of the original $I_d$: the splitting planes are still the same seven labelled planes, whether or not new quadratic splitting subfields appear. This uses the functoriality of the comparisons and the equality of the kernels on $\A$ and $\Lam$ from Section~\ref{sec:surface}. Induction on the number of factors gives the claimed kernel on the original norm quotient: at each step divide a class killed by restriction by a lift in $I_d$. The comparisons between the seven planes still generate the kernel for the conic, so restriction from $J_C$ onto it is surjective. In dimension six the calculation with the augmentation space in the preceding proof gives the dimension of the sum as the number of planes meeting $D$ minus one.
\end{proof}

\section{A rank criterion over the ground field}\label{sec:rank}
Relations between the comparison classes can be tested by quaternion symbols over any quadratic subfield. For any quadratic $K/F$ in $L$, write $L=K(\sqrt r,\sqrt s)$ with $r,s\in F^\times$ and put
\[
 \begin{split}
 \Phi_K(v)&=((r,v)_K,(s,v)_K),\quad v\in K^\times,\\
 D_K&=\Phi_K(F^\times),\\
 W_K&=\{(\alpha,\beta)\in\Phi_K(K^\times):
 \cor_{K/F}\alpha=\cor_{K/F}\beta=0\}.
 \end{split}
\]
These are elementary abelian $2$-groups, and $D_K\subset W_K$ by the projection formula.
\begin{proposition}\label{prop:rank}
There is an isomorphism
\[
 W_K/D_K\xrightarrow{\ \sim\ }\Lam(L/F),\qquad
 \Phi_K(v)\longmapsto [\N_{K/F}(v)].
\]
Let $c_1,\ldots,c_n\in F^\times$ represent the generators of a kernel $\K(L/F;F(q))$ obtained in Sections~\ref{sec:surface} and~\ref{sec:conic}. Choose $v_i\in K^\times$ with $\N_{K/F}(v_i)=c_i$ and set $U_K=\langle\Phi_K(v_i):1\le i\le n\rangle$. Then
\begin{equation}\label{eq:rank}
 \dim_{\F_2}\K(L/F;F(q))=\dim_{\F_2}U_K-\dim_{\F_2}(U_K\cap D_K).
\end{equation}
The right side is independent of the lifts and of the choice of $K,r,s$.
\end{proposition}
\begin{proof}
For the biquadratic norm identity (see \cite[Lemma 3.3]{AB}), suppose $u$ is a norm from both $K(\sqrt r)$ and $K(\sqrt s)$ and choose $x$ and $y$ in these respective fields with norms $u$. With bars denoting quadratic conjugation,
\[
 \N_{L/K}(x+y)=u(x+\bar x+y+\bar y)^2.
\]
For the identity, multiply first $(x+y)(x+\bar y)=x(x+\bar x+y+\bar y)$ and then its conjugate in $K(\sqrt r)$. Replacing $y$ by $-y$ if necessary gives a nonzero trace sum: both traces cannot vanish, as then $u=-r x_0^2=-s y_0^2$ for nonzero $x_0,y_0\in K$, contradicting the independence of $r,s$ over $K$. Hence $u\in K^{\times2}\N_{L/K}(L^\times)$, and the reverse inclusion follows by norm transitivity. Therefore
\[
 \ker\Phi_K=K^{\times2}\N_{L/K}(L^\times).
\]
Thus $\Phi_K$ identifies $\A(L/K)$ with its image. The norm induces a homomorphism $\A(L/K)\to\A(L/F)$ whose kernel is the image of $F^\times$. Indeed, if $\N_{K/F}(v)=c^2\N_{L/F}(w)$, with $c\in F^\times$ and $w\in L^\times$, then $v/(c\N_{L/K}(w))$ has norm one. By Hilbert 90 it equals $h/\bar h$ for some $h\in K^\times$, and modulo squares this equals $h\bar h\in F^\times$. The converse is immediate since a constant from $F$ has square norm.

By the projection formula, $[\N_{K/F}(v)]$ belongs to $\Lam(L/F)$ precisely when $\Phi_K(v)$ has both corestrictions zero. Every representative of a class of $\Lam(L/F)$ is a norm from $K$, so the induced map on $W_K/D_K$ is surjective and hence an isomorphism. The rank formula follows by restricting this map to $U_K$. The dimension is intrinsic because changing a lift changes its image only by an element of $D_K$, and the isomorphism holds for every $K,r,s$.
\end{proof}
In particular, for nonsquare discriminant with $K=F(\sqrt{\disc q})\subset L$ and all three quadratic subfields of $L/K$ splitting the division algebra $Q=C_0(q)$, one can take $U_K=\langle Q\rangle\oplus\langle Q\rangle$. Indeed common slots for $(r,s)$ and $(r,rs)$ have symbol pairs $(Q,Q)$ and $(Q,0)$ respectively. Thus the kernel is two-dimensional precisely when no nonzero pair in $U_K$ comes from the same element of $F^\times$. More generally \eqref{eq:rank} reduces a relation among comparison classes to a pair of symbol equations with one unknown $c\in F^\times$. For square discriminant the representatives $c_i$ are obtained from \eqref{eq:pairnorm}, and one may use any fixed quadratic $K\subset L$ for their lifts.

In particular, if $Q=(r,u)_K=(s,u)_K$ for a comparison between $K(\sqrt r)$ and $K(\sqrt s)$, its class is trivial if and only if there is a common slot defined over $F$: that is, there is $c\in F^\times$ with $Q=(r,c)_K=(s,c)_K$. This follows by applying the isomorphism to $\Phi_K(u)=(Q,Q)$. If all three quadratic extensions of $K$ split $Q$, the product of the three pairwise comparisons is $1$. The dimension of the kernel is two, one, or zero according as none, exactly one, or all three of the conditions for descent of a common slot hold.

\section{Sharpness}\label{sec:sharp}
In this section write $k_{u,v}=k(\sqrt u,\sqrt v)$ for a biquadratic extension of a base field $k$.

\begin{proposition}\label{prop:onehot}
Let $k_*=\mathbb Q(\sqrt{-3})(x)$, let $a=x+4$, $b=x+1$, $e=x$, and put $L_*=k_*(\sqrt a,\sqrt b,\sqrt e)$. The quaternion algebra
\[
 Q=(ab,3)=(a,-b)
\]
is division and has the following properties.
\begin{enumerate}
\item All seven biquadratic subfields split $Q$. The quadratic splitting classes are exactly $a,ab,e,be$, the complement of the plane spanned by the classes of $b,ae$.
\item Relative to the reference $k_*(\sqrt a)$, the comparison $c_{M,k_*(\sqrt a)}$ is $\xi=[-x]\ne1$ if $M=(k_*)_{b,ae}$, and is trivial for every other biquadratic subfield $M$.
\end{enumerate}
\end{proposition}
\begin{proof}
The radicands are independent by the valuations at their three distinct zeros. We first prove the nonvanishing of $\xi$ in $\A(L_*/k_*)$, using the dyadic argument of Shapiro, Tignol and Wadsworth \cite[Theorem 5.1 and Remark 5.4]{STW}. This verifies that the class of Gille \cite[Remark 4]{Gille} survives the extension of the field of constants.

Put $k_0=\mathbb Q(\sqrt{-3})$, $K=k_*(s)$, with $s^2=x$. The prime $2$ is inert and unramified in $k_0$ (the polynomial $T^2+T+1$ is irreducible modulo $2$). Equip $K=k_0(s)$ with the Gauss extension of its dyadic valuation, normalised by $v(2)=1$. Its residue field is $\F_4(\bar s)$; the residue field of the restricted valuation on $k_*$ is $\F_4(\bar s^2)$. We have
\[
 (s+2)^2-a=4s,\qquad (s+1)^2-b=2s.
\]
If $[-x]=[\N_{K/k_*}(s)]$ were trivial in $\A(L_*/k_*)$, the description of the kernel of the norm in Proposition~\ref{prop:rank} would give $f\in k_*^\times$ with $s/f\in K^{\times2}\N_{L_*/K}(L_*^\times)$. Such an element is a norm from both quadratic subfields over $K$. The two identities above show that $s$ is a norm from $K(\sqrt a)$ and $2s$ is a norm from $K(\sqrt b)$, so $f$ and $2f$, respectively, are norms from these fields.

Let $\widehat K$ be the completion of $K$ for $v$. The element $(\sqrt a/s-1)/2$ satisfies $T^2+T=1/s^2$. Its reduction is irreducible over $\F_4(\bar s)$: modulo the Artin--Schreier image, $1/\bar s^2$ is equivalent to $1/\bar s$, which has a simple pole. Thus $\widehat K(\sqrt a)/\widehat K$ is unramified of degree two, and the norm condition on $f$ implies $v(f)\in2\Z$.

By contrast, $\widehat K(\sqrt b)/\widehat K$ is totally ramified of degree two: $\delta=s+1-\sqrt b$ satisfies the Eisenstein polynomial
\[
 T^2-2(s+1)T+2s.
\]
Let $N_b$ denote the norm from $\widehat K(\sqrt b)$ to $\widehat K$. Then $N_b(\delta)=2s$. For the extended valuation $w$ normalised by $w(\delta)=1$, we have $v(N_b(y))=w(y)$. Suppose $c=N_b(y)\in k_*^\times$ and put $n=v(c)$. Then $c/(2s)^n=N_b(y/\delta^n)$ is the norm of a unit. Its residue is a square in the unchanged residue field $\F_4(\bar s)$, and equals
\[
 \overline{c/2^n}\,\bar s^{-n},\qquad
 \overline{c/2^n}\in\F_4(\bar s^2)^\times=\F_4(\bar s)^{\times2}.
\]
Since $\bar s$ is not a square, $n$ is even. This parity assertion concerns norms in $k_*^\times$, not all norms from the ramified extension. Applying it to $c=2f$ contradicts $v(2f)=1+v(f)$. Hence $\xi\ne1$ in $\A(L_*/k_*)$.

The Steinberg relation $(a/3,-b/3)=0$ and $(3,-3)=0$ give $(a,-b)=(a,3)+(b,3)=(ab,3)$. At the valuation $x=-4$, the tame residue (see \cite[Chapter 6]{GS}) of $(ab,3)$ is the class of $3$, which is nonsquare in $k_0$. Hence $Q$ is division. Thus $a$ and $ab$ are splitting classes. The class $be$ splits $Q$ since over $k_*(\sqrt{be})$ one has $(a,-b)=(a,-e)=0$, using $a=4+e$. Also $Q=(ab,-1)$ because $-3$ is a square in $k_*$. Over $k_*(\sqrt e)$, both $a$ and $b$ are sums of two squares, so $e$ is a splitting class. The four splitting classes so obtained form the complement of the plane spanned by the classes of $b,ae$. Every other plane contains one of them.

For the remaining plane let $K'=k_*(\sqrt{ae})$ and put $h=\sqrt{ae}/e$, $D_-=h^2-1$, and $D_+=h^2+3$. We have $x=4/D_-$, $a=4h^2/D_-$ and $b=D_+/D_-$. In particular the classes of $a,b$ over $K'$ are those of $D_-,D_-D_+$, respectively, and $Q_{K'}=(D_-,D_+)$. Put $u=2D_+(1-h)$. Then
\[
 u=D_+^2-D_+(h+1)^2,\qquad u/D_+=(h-1)^2-D_-.
\]
These identities give $(D_+,u)=0$ and $(D_-,u)=(D_-,D_+)$, whence
\[
 Q_{K'}=(a,u)_{K'}=(b,u)_{K'}.
\]
Thus $(k_*)_{b,ae}$ splits $Q$. Its comparison with $(k_*)_{a,e}$, using their intersection $K'$, is given by \eqref{eq:pairnorm}: indeed the two fields are $K'(\sqrt b)$ and $K'(\sqrt a)$. Since $h$ is sent to $-h$ by conjugation,
\[
 \N_{K'/k_*}(u)=-4D_+^2D_-=-16D_+^2/x.
\]
Here $D_+=h^2+3$ belongs to $k_*$. The norm has square class $-x$. The comparison of $(k_*)_{a,e}$ with $k_*(\sqrt a)$ is trivial, so the exceptional comparison equals $\xi$ and belongs to $\Lam(L_*/k_*)$. All other comparisons are trivial, since the corresponding fields contain a quadratic splitting subfield. If a quadratic subfield of the remaining plane split $Q$, its comparison both with that plane and with $k_*(\sqrt a)$ would be trivial, forcing $\xi=1$. This proves that the four listed quadratic splitting classes are the only ones.
\end{proof}

\begin{lemma}\label{lem:specialize}
Let $\Gamma$ be a smooth geometrically integral curve over $E$, with $P\in\Gamma(E)$ and $\operatorname{char} E\ne2$. The embedding of $E(\Gamma)$ into the completion $E((z))$ at $P$ identifies the square class of every function regular and nonzero at $P$ with the class of its value at $P$. For any finite Galois extension $M/E$, the natural map
\[
 \A(M/E)\longrightarrow\A(M((z))/E((z)))
\]
is injective.
\end{lemma}
\begin{proof}
The completed local ring at $P$ is $E[[z]]$. If $v$ is a unit of the local ring, then $v/v(P)$ belongs to $1+zE[[z]]$ and has a square root. For the second assertion, taking leading coefficients (denoted by $\lc$) in an equality $c=h^2\N_{M((z))/E((z))}(f)$ with $c\in E^\times$ gives $c=\lc(h)^2\N_{M/E}(\lc(f))$.
\end{proof}

For the rest of this section let $V=\F_2^3$ and take $k_*$ and $Q$ from Proposition~\ref{prop:onehot}, with $C=\SB(Q)$.

\begin{lemma}\label{lem:local-seed}
There is an embedding $k_*\hookrightarrow\mathbb Q_7$ for which $Q_{\mathbb Q_7}$ is division.
\end{lemma}
\begin{proof}
The element $-3$ is a square in $\mathbb Q_7$ by Hensel's lemma. Choose the image of $x$ in the nonempty open set $v_7(x+4)=1$ to be transcendental over $\mathbb Q$ (only countably many elements are algebraic over $\mathbb Q$). This defines the embedding of $k_*$. Then $x+1$ has residue $4$, so $ab$ has valuation one. Since $3$ is a nonsquare unit, $Q=(ab,3)$ remains division over $\mathbb Q_7$.
\end{proof}

\begin{lemma}\label{lem:parameter}
Let $S_0=(\mathbf G_m)^3$ over $k_*$ with coordinates $t_1,t_2,t_3$. For each nonzero subspace $W\subseteq V$ let $S_W\to S_0$ be the finite \'etale cover obtained by adjoining square roots of the monomials for a basis of $W$. For any collection $\mathcal M$ of such subspaces put
\[
 Y=\prod_{W\in\mathcal M}\Res_{S_W/S_0}(C\times_{k_*}S_W),
\]
with product over $S_0$. The variety $Y$ is smooth, quasi-projective and geometrically integral over $k_*$. Any finite set of closed points of $Y$ is contained in a smooth geometrically integral curve $\Gamma\subset Y$ defined over $k_*$. At a marked point $P$ of this curve, $k_*(\Gamma)$ embeds in $k_*(P)((z))$ so that the square classes of the radicands are the classes of their values at $P$.
\end{lemma}
\begin{proof}
The Weil restrictions are smooth and projective over the torus, with geometric fibres products of projective lines. Hence $Y$ is smooth and geometrically integral (for the empty collection, $Y=S_0$).

Let $Z$ be the reduced union of the marked closed points. It is geometrically reduced since $k_*$ has characteristic zero. In a locally closed embedding of $Y$ in projective space, use the linear system of hypersurfaces of degree $m$ containing $Z$. For $m$ sufficiently large it separates points and tangent directions off $Z$, and its first derivatives span the cotangent space at every geometric point of $Z$ (by separation of jets for sufficiently high powers of the hyperplane bundle). If $\dim Y\ge2$, the map off the base locus therefore has image of dimension at least two. The Bertini smoothness and irreducibility theorems give a smooth geometrically integral section off $Z$; requiring a nonzero derivative at each point of $Z$ ensures smoothness there as well. No additional component can be supported on the finite set $Z$, since the section has pure codimension one. These are nonempty open conditions on the projective space of hypersurfaces defined over $k_*$ (see \cite{Bertini}); since $k_*$ is infinite, one can choose a hypersurface over $k_*$. Repeating on the resulting section gives $\Gamma$.

At a marked point $P$, write $E_P=k_*(P)$. Hensel's lemma lifts a primitive element of $E_P/k_*$ to the completed discrete valuation ring, giving a coefficient field over $k_*$. A uniformiser $z$ identifies this ring with $E_P[[z]]$. Each $t_j$ is a unit and $t_j/t_j(P)\in1+zE_P[[z]]$ is a square. This gives the asserted embedding and square classes.
\end{proof}

\begin{lemma}\label{lem:local-point}
Let $k'$ be a dense subfield of a nonarchimedean local field $H$, and let $Z$ be a smooth variety over $k'$. In every nonempty open subset of $Z(H)$ for the local field topology there is an $H$-point whose image on $Z$ is a closed point.
\end{lemma}
\begin{proof}
Choose an $H$-point in the given open set and an \'etale chart $U\to\mathbb A^n_{k'}$ around its image on $Z$. By the local inverse function theorem there is a smaller neighbourhood of the $H$-point, contained in the given open set, mapping homeomorphically onto an open subset of $H^n$. By density this subset contains a point of $k'^n$. Its lift lies in a fibre of a finite type \'etale morphism over a $k'$-rational point. The fibre is zero-dimensional, so the lift factors through a closed point with residue field finite over $k'$ and embedded in $H$.
\end{proof}

\begin{lemma}\label{lem:number-pattern}
Every collection $\mathcal S$ of nonzero subspaces of $V=\F_2^3$ closed under passing to larger subspaces occurs as the splitting pattern of a quaternion division algebra in a triquadratic extension of a number field. In particular, if $b(S,\mathcal T)=0$, the bound in \eqref{eq:pattern-bound} is attained over a number field.
\end{lemma}
\begin{proof}
Declare the zero subspace to be omitted and list the maximal omitted subspaces $U_1,\ldots,U_m$. Take $k=\mathbb Q(\sqrt{-7})$. For each $U_i\ne0$ choose a pair of distinct odd finite places of $k$, all pairs disjoint. At each place $v$ in this pair choose a linear map
\[
 \psi_v:V\longrightarrow k_v^\times/k_v^{\times2}
\]
with kernel $U_i$, possible since the local space of square classes has dimension two. If $U_i=0$, it is the only maximal omitted subspace. Since $-7\equiv1\pmod8$, the prime $2$ splits in $k$ and both completions are $\mathbb Q_2$. Use these two places and choose injective maps into their three-dimensional spaces of square classes. If no dyadic place was used, choose an injective map at one of the two places above $2$ as well, to ensure independence of the radicands. This auxiliary place is not in the ramification set of $Q$.

By weak approximation choose three elements of $k^\times$ representing the specified images of a fixed basis at all these places; the conditions on square classes are open. They are independent at a dyadic place, hence over $k$, and define a triquadratic extension $L/k$ with space of square classes $V$. The specified pairs give a nonempty even set of places. The Brauer--Hasse--Noether theorem gives a class with invariant $1/2$ at these places and zero elsewhere. Its index is two (the index over a number field is the least common multiple of the local indices), so its division representative is a quaternion algebra $Q/k$ (see \cite{GS}). For each $W\subseteq V$, write $k_W$ for the corresponding subfield of $L/k$. At a place of $k_W$ above $v$ the local extension has degree $2^{\dim\psi_v(W)}$. Therefore the invariant $1/2$ is killed at every place over $v$ precisely if $W\not\subseteq U_i$. Consequently $Q_{k_W}$ splits precisely if $W$ is contained in no maximal omitted subspace, that is, precisely if $W\in\mathcal S$. For number fields $\Lam(L/k)=0$ by \cite[Theorem 3.7]{STW}, so $J_C=0$.
\end{proof}

\begin{proof}[Proof of Theorem~\ref{thm:pattern}]
If $S$ is empty, the comparisons between the $|\mathcal T|$ splitting planes generate the kernel, with at most $|\mathcal T|-1$ generators when $\mathcal T$ is nonempty. If $S$ is nonempty, all quadratic splitting fields have trivial comparisons with one another, as do all splitting planes containing one of them. Thus only the $|\mathcal T_0|$ comparisons with a fixed quadratic splitting field are needed, giving the upper bound. The case $b(S,\mathcal T)=0$ follows from Lemma~\ref{lem:number-pattern}. Assume henceforth $b(S,\mathcal T)>0$. Take $\mathcal M$ to be the set of minimal members of $\mathcal S$ in Lemma~\ref{lem:parameter}. Let $S$ and $\mathcal T$ be the sets prescribed by $\mathcal S$; in particular all planes containing a member of $S$ belong to $\mathcal T$.

If $S$ is empty, fix a reference plane $W_0\in\mathcal T$. For each $W_i\in\mathcal T\setminus\{W_0\}$ choose an isomorphism from $V$ onto the space of square classes of $L_*/k_*$ sending $W_i$ to the exceptional plane. Its values on a basis determine a point of $S_0(k_*)$. All required splitting fields split $Q$ there, so this point lifts to $Y(k_*)$. At this point the comparison $c_{W_i,W_0}$ specialises to $\xi$ and all the other comparisons with $W_0$ specialise to $1$. If $S$ is nonempty, fix $d_0\in S$. For each $W_i\in\mathcal T_0$ choose an isomorphism sending $W_i$ to the exceptional plane. Every element of $S$ then maps outside that plane and hence is a quadratic splitting class of the example. Again the point lifts to $Y(k_*)$, and the comparisons of $W_i$ with the reference quadratic field have values $\xi$ for the indexed plane and $1$ for all others. These marked points are distinct because their isomorphisms send different planes to the exceptional plane.

For the omitted subspaces we impose local conditions. Embed $k_*$ into $H=\mathbb Q_7$ as in Lemma~\ref{lem:local-seed}; its image is dense since it contains $\mathbb Q$. For each nonzero $U\notin\mathcal S$, choose a linear map
\[
 V\longrightarrow H^\times/H^{\times2}
\]
with kernel $U$. Call the map $\psi_U$; it exists since $\dim(V/U)\le2=\dim_{\F_2}(H^\times/H^{\times2})$. Choose values of the $t_j$ representing the images of a basis. For $W\in\mathcal M$, the finite \'etale algebra of specialised roots for $W$ is a product of $2^{\dim W-\dim\psi_U(W)}$ copies of a field of degree $2^{\dim\psi_U(W)}$ over $H$. Since $W\not\subseteq U$, this degree is even, so each factor splits $Q_H$ by the formula for restriction of local Brauer invariants. The defining property of Weil restriction therefore gives a lift of the base point to $Y(H)$. The square classes of the radicands are constant in a sufficiently small local neighbourhood of the lift. Lemma~\ref{lem:local-point} gives a closed point in this neighbourhood with residue field embedded in $H$ and with the same map $\psi_U$ on square classes.

Take a curve $\Gamma$ through all the marked points by Lemma~\ref{lem:parameter}, and set $F=k_*(\Gamma)$ and $L=F(\sqrt{t_1},\sqrt{t_2},\sqrt{t_3})$. There is at least one marked point from the example, since $b(S,\mathcal T)>0$, so independence of the radicands follows by completion at that point. Thus $L/F$ is triquadratic. At the generic point of the curve, the defining property of Weil restriction gives a point of $C$ over each $F_W$ with $W\in\mathcal M$ (these \'etale algebras are now fields by independence). Thus $Q$ splits over every prescribed subfield.

For each marked point $P$ for an omitted $U$, Lemma~\ref{lem:parameter} gives $F\hookrightarrow k_*(P)((z))\hookrightarrow H((z))$, the second map extending $k_*(P)\hookrightarrow H$. For every vector in $U=\ker\psi_U$, the corresponding monomial in the $t_j$ has square residue in $H$ and is thus a square in $H((z))$. The embedding of $F$ extends to $F_U$. Since $Q_H$ remains division over $H((z))$ by injectivity of the constant Brauer group, $Q_{F_U}$ is nonsplit. If every nonzero subspace were prescribed to split $Q$, then $\mathcal T_0$ would be empty and $b(S,\mathcal T)=0$. Thus there is at least one omitted nonzero subspace, and in particular $Q_F$ is division.

At each marked $k_*$-point from the example, Lemma~\ref{lem:specialize} identifies $L\otimes_F k_*((z))$ with the corresponding constant extension $L_*((z))$, with the subfields labelled by the chosen isomorphism of spaces of square classes. By functoriality of the comparisons, the generator indexed by this point restricts to $\xi$, and all other generators restrict to $1$, in $\A(L_*((z))/k_*((z)))$. The nonzero class $\xi$ remains nonzero there by Lemma~\ref{lem:specialize}. Applying these restrictions to any relation proves independence of the $b(S,\mathcal T)$ generators. This gives equality in \eqref{eq:pattern-bound}, and $\operatorname{tr.deg}_{\mathbb Q}F=2$.
\end{proof}

\begin{proof}[Proof of Corollaries~\ref{thm:square-intro} and~\ref{thm:sharp-intro}]
Prescribe all seven splitting planes and no splitting lines in Theorem~\ref{thm:pattern}. Then $\dim_{\F_2} J_C=6$ and Theorem~\ref{thm:compatibility} gives $\dim_{\F_2} I_d=2$ for every $d\ne0$; each $q_d$ is anisotropic since $Q_{F(\sqrt{\widetilde d})}$ is division. The field has transcendence degree two. For the other positive bounds prescribe all seven planes and $r$ independent splitting classes, $r=1,2,3$. Exactly $2^{3-r}$ planes avoid those classes, so Theorem~\ref{thm:pattern} attains the bound. For $r=3$ the example of Proposition~\ref{prop:onehot} already suffices.
\end{proof}

\begin{remark}
Every dimension $0\le j\le6$ occurs over a field of transcendence degree at most two over $\mathbb Q$: for $j>0$ prescribe exactly $j+1$ planes, no lines and the full space; for $j=0$ take a number field. Proposition~\ref{prop:onehot} gives a nonzero kernel in transcendence degree one. The least transcendence degree for the largest kernel remains undetermined.
\end{remark}

\section*{Declarations}
\noindent\textbf{Funding.} This research was supported by the National Natural Science Foundation of China (NSFC), Grant Nos.~12231009 and~11971224.

\medskip\noindent\textbf{Conflict of interest.} The author declares that he has no conflict of interest.

\medskip\noindent\textbf{Data availability.} No research datasets were generated or analysed for this study.


\begin{thebibliography}{99}
\bibitem{AB} M. Archita and K. J. Becher, Rational connectedness for groups of proper projective similitudes, arXiv:2506.21717 (2025).
\bibitem{BCS} M. Borovoi, J.-L. Colliot-Th\'el\`ene and A. N. Skorobogatov, The elementary obstruction and homogeneous spaces, \emph{Duke Math. J.} \textbf{141} (2008), 321--364.
\bibitem{CTS} J.-L. Colliot-Th\'el\`ene and A. N. Skorobogatov, \emph{The Brauer--Grothendieck Group}, Ergebnisse der Mathematik und ihrer Grenzgebiete, vol. 71, Springer, Cham, 2021.
\bibitem{EKM} R. Elman, N. Karpenko and A. Merkurjev, \emph{The Algebraic and Geometric Theory of Quadratic Forms}, American Mathematical Society Colloquium Publications, vol. 56, American Mathematical Society, Providence, RI, 2008.
\bibitem{Gille} P. Gille, Examples of non-rational varieties of adjoint groups, \emph{J. Algebra} \textbf{193} (1997), 728--747.
\bibitem{GS} P. Gille and T. Szamuely, \emph{Central Simple Algebras and Galois Cohomology}, second edition, Cambridge Studies in Advanced Mathematics, vol. 165, Cambridge University Press, Cambridge, 2017.
\bibitem{Bertini} J.-P. Jouanolou, \emph{Th\'eor\`emes de Bertini et applications}, Progress in Mathematics, vol. 42, Birkh\"auser, Boston, 1983.
\bibitem{LW} D. B. Leep and A. R. Wadsworth, The transfer ideal of quadratic forms and a Hasse norm theorem mod squares, \emph{Trans. Amer. Math. Soc.} \textbf{315} (1989), 415--431.
\bibitem{Siv18} A. S. Sivatski, On the homology groups of the Brauer complex for a triquadratic field extension, \emph{Math. Nachr.} \textbf{291} (2018), 518--538.
\bibitem{Siv27} A. S. Sivatski, Quadratic forms isotropic over two quadratic extensions, \emph{J. Algebra} \textbf{715} (2027), 303--321. \href{https://doi.org/10.1016/j.jalgebra.2026.09.043}{doi:10.1016/j.jalgebra.2026.09.043}.
\bibitem{STW} D. Shapiro, J.-P. Tignol and A. R. Wadsworth, Witt rings and Brauer groups under multiquadratic extensions. II, \emph{J. Algebra} \textbf{78} (1982), 58--90.
\bibitem{Tignol} J.-P. Tignol, Corps \`a involution neutralis\'es par une extension ab\'elienne \'el\'ementaire, in \emph{Groupe de Brauer} (Les Plans-sur-Bex, 1980), Lecture Notes in Mathematics, vol. 844, Springer, Berlin, 1981, pp. 1--34.
\end{thebibliography}
\end{document}